\documentclass{amsart} 
\theoremstyle{plain}
\newtheorem{thrm}{Theorem}[section]
\newtheorem{lmm}[thrm]{Lemma}

\newtheorem{prpstn}[thrm]{Proposition}
\theoremstyle{definition}
\newtheorem{dfntn}[thrm]{Definition}
\newtheorem{rmrk}[thrm]{Remark}
\theoremstyle{plain}
\newcommand{\ie}{\emph{i.e.}}
\newcommand{\cf}{cf.}
\def\xR{\mathbb{R}}
\def\xZ{\mathbb{Z}} 
\def\xCtwo{{\rm C}^{2}} 
\usepackage[latin1]{inputenc}
\usepackage[T1]{fontenc}
\usepackage{url}
\newcommand{\newspin}[2]{%
  \expandafter\newcommand\csname spin#1\endcsname[1]{#2_{##1}}%
  \expandafter\newcommand\csname spins#1\endcsname{\mathbf{#2}}%
  \expandafter\newcommand\csname spins#1Boite\endcsname[1]{\mathbf{#2}_{##1}}%
  \expandafter\newcommand\csname espL#1\endcsname{\esp^{\Lambda,\mathbf{#2}}}%
}
\newspin{}{x}
\newspin{Y}{y}
\newspin{Z}{z}

\newcommand*{\term}[1]{\emph{#1}}

\newcommand*\Ent{\mathbf{Ent}}
\newcommand*\Var{\mathbf{Var}}

\newcommand*{\abs}[1]{\left|#1\right|}
\newcommand*{\norme}[2]{\left\|#2\right\|_{#1}}
\newcommand*{\ninf}[1]{\norme{\infty}{#1}}
\newcommand*{\alphaNorm}[1]{\norme{\alpha}{#1}}

\newcommand*{\lDeuxAlpha}{\ell^2(\alpha)} 
\newcommand*{\temperedDLRStates}{\mathcal{G}_t} 
\newcommand*{\wellTemperedCfgs}{\mathcal{P}} 
\newcommand*{\convol}{\star} 
\newcommand*{\esp}{\mathbf{E}}
\newcommand*{\ind}[1]{\mathbf{1}_{#1}}
\DeclareMathOperator{\card}{Card}
\DeclareMathOperator{\osc}{osc}
\newcommand*{\egalDef}{\overset{\text{\tiny\textit{def}}}{=}}

\newcommand\SIS{\ensuremath \mathcal{S}'}
\newcommand\temperedCfgs{\SIS}

\newcommand{\relEnt}[2]{I(#1 | #2)}
\newcommand{\loi}[1]{\mathcal{L}(#1)}
\newcommand*{\GBIa}{\textbf{GBI}($a$)}
\newcommand*{\GBIaPrime}{\textbf{GBI}($a'$)}

\newcommand{\xBoite}{X_t^{L,0,\spins}}
\newcommand{\xInf}{X_t^{\spins}}

\newcommand{\cNt}{c^{LS}_{n(t)}}
\newcommand{\dCfgs}{d_{tmp}}
\newtheorem{hpthss}{Hypothesis}
\newtheorem{assmptn}{Assumption}

\newcommand*{\eg}{\emph{e.g.}}

\begin{document}
\title[Functional inequalities and uniqueness of the Gibbs measure]{Functional inequalities and uniqueness of the Gibbs measure --- from log-Sobolev to Poincaré}%
\author{Pierre-André Zitt}
\address{Équipe Modal'X, EA3454\\ Université Paris X, Bât. G\\ 200 av. de la République\\ 92001 Nanterre\\ France}
\email{pzitt@u-paris10.fr}
\begin{abstract}
  In a statistical mechanics model with unbounded spins, we prove uniqueness of the Gibbs measure
under various assumptions on finite volume functional inequalities. We follow Royer's approach (\cite{Roy99})
and obtain uniqueness by showing convergence properties of a Glauber-Langevin dynamics. The result was known when 
the measures on the box $[-n,n]^d$ (with free boundary conditions) satisfied the same logarithmic Sobolev inequality. We generalize this in two directions: either the constants may be allowed to grow sub-linearly in the diameter, or we may suppose a weaker inequality than log-Sobolev, but stronger than Poincaré. We conclude by giving a heuristic argument showing that this could be the right inequalities to look at. 
\end{abstract}
%
%
\maketitle
\section*{Introduction}
Questions of convergence of dynamical models of statistical physics (\eg{} Glauber dynamics for the classical Ising model) have prompted people to study \emph{functional inequalities} for the equilibrium measures related to these dynamics, \ie{} for Gibbs states. These inequalities are indeed a good way to obtain convergence results for semi-groups. 
Moreover, if the classical functional inequalities (Poincaré, logarithmic Sobolev) are known to \emph{tensorize} in a good way, studying them for non-product measures in large dimensions was much more challenging, and Gibbs measures are a natural example of these non-product measures. Therefore, many authors (see \eg{} \cite{LY93,SZ92c} for the bounded spins case, \cite{BH99,Yos01,Led01} for the unbounded case) have investigated links between ``uniform'' functional inequalities, convergence of associated dynamics and mixing properties of equilibrium measures.

In several cases, it was also proved that there is a regime in which all these ``good'' properties hold simultaneously. 

The ``uniformness'' we alluded to is typically ``uniform on all (regular) finite sets, and all boundary conditions''. However, in his book \cite{Roy99}, G.~Royer shows that a logarithmic Sobolev inequality, uniform over the boxes $[-n,n]^d$, for a \emph{single} boundary condition entails the uniqueness of the infinite volume Gibbs measure. 

We show here that this assumption may be relaxed in two different ways. Firstly, we show that the constants may be allowed to grow sublinearly in $n$. Secondly, we may replace logarithmic Sobolev inequalities by weaker inequalities, and show the uniqueness when we only suppose uniform Beckner inequalities (\cf{} theorem \ref{thm=uniciteParBeckner} for a precise statement). 

After introducing notations and preliminary estimates (section \ref{sec=intro}), we prove the result concerning logarithmic Sobolev inequalities (sec. \ref{sec=LS}). In the last section, we show the result on Beckner inequalities and indicate a heuristic argument that they may be the critical scale for uniqueness.

\section{Notations and preliminary estimates}
\label{sec=intro}
\subsection{The model~: equilibrium and dynamics}
\subsubsection{The model --- equilibrium}
We consider a variant of the classical Ising model. To define it, we briefly introduce the following notions, referring to \cite{Zit06} for details. 
A \term{configuration} is a map $\spins: \mathbb{Z}^d \to \mathbb{R}$. We denote by $\spins_L$ the restriction of $\spins$ to the subset $L\subset \mathbb{Z}^d$. When $L$ is a singleton, we will simply write $\spin{i}$; this is the \term{spin} at site $i$. To each finite subset $L$ of $\mathbb{Z}^d$, and each configuration $\spinsZ$ (\term{boundary condition}), we associate a \term{Hamiltonian} 
\[
U_{L,\spinsZ}(\spins) = \sum_{i\in L} V(\spin{i}) + \sum_{i,j\in L, i\sim j} J(i-j) \spin{i} \spin{j}
+ \sum_{i\in L, j\notin L} J(i-j) \spin{i}\spinZ{j}.
\]
where $V$ and $J$ satisfy the following.
\begin{hpthss}[Self-interaction]
  \label{hpth=formeDeV}
  The function $V$ satisfies:
  \begin{itemize}
    \item convexity at infinity --- there exists $V_1$, $V_2$ such that $V=V_1 = V_2$,  $V_2$ is $\xCtwo$ and compactly supported, $\inf V_1'' >0$. 
    \item Polynomial growth --- there exists constants $a_V,b_V$, and a  $d_V>0$ such that for all $x$, $\abs{V(x)}\leq a_V \abs{x}^{d_V} + b'_V$.
    \item There exists a  $ a'_V\in]0,1[ $ such that  $x\mapsto a'_V V'(x)^2 -V''(x) $ is bounded from below.
  \end{itemize}
\end{hpthss}
\begin{hpthss}
  \label{hpth=formeDeJ}
  The interaction $J:\mathbb{Z}^d \to \xR$ is a symmetric function with finite support. We also define $p(i) = \abs{J(i)}$ 
  and suppose that 
  \[ \sigma \egalDef \sum_{i\in\xZ^d} p(i) < \inf V_1'',\]
  where $V_1$ is defined by the previous hypothesis.
\end{hpthss}
These hypotheses are satisfied for the usual models, namely the gaussian case and the 
double well potential~:
\[
  V(x) = ax^4 - bx^2.
\]

We then define the \term{finite volume Gibbs measure} on $\mathbb{R}^L$ by~:
\[
  d\mu_{L,\spinsZ}(d\spins_L) = \frac{1}{Z_{L,\spinsZ}} \exp\left( - U_{L,\spinsZ}(\spins_L) \right).
\]
where $Z_{L,\spinsZ} = \int\exp(-U_{L,\spins})d\spins_L$ 
is a normalizing constant (note that we may abuse notations and speak of $U_{L,\spinsZ}(\spins_L)$, since $U_{L,\spinsZ}$ only depends on the spins in $L$). 

Note that $\mu_{L,\spinsZ}(d\spins)$, while it is originally defined as a measure on $\mathbb{R}^L$, may be extended to $\mathbb{R}^{\mathbb{Z}^d}$ by fixing $\spins=\spinsZ$ outside $L$. This enables us to define a kernel $\mu_L$~:
\[
  \mu_L : \spinsZ \mapsto \mu_{L,\spinsZ}.
\]
An \term{infinite volume} Gibbs measure is a measure on $\mathbb{R}^{\mathbb{Z}^d}$ that satisfies the \term{DLR equations}~:
\[
  \forall L, L \text{ finite}, \quad \mu\mu_L = \mu.
\]
For technical and physical reasons, we will only consider \term{tempered} configurations and measures. 
For every $d$, let $\wellTemperedCfgs(d)$ be defined by
\[ 
  \spins \in \wellTemperedCfgs(d)
  \iff \exists c_\spins, \forall i, \quad \spin{i} \leq (1+\abs{i})^d.
\]
A configuration is called \term{tempered} if it is in $\wellTemperedCfgs(d)$, for some $d$.
A tempered measure is one that satisfies~:
\[
 \exists C_\mu, \forall i\in\mathbb{Z}^d, \quad
   \mu(\abs{\spin{i}} )\leq C_\mu.
\]
It may be shown (\cite{Zit06}, section 2.2) that this is equivalent to other standard definitions of temperedness, and that there exists a $\dCfgs$, depending only on the dimension $d$, such that every tempered measure $\mu$ satisfies:
\begin{equation}
  \mu\left( \wellTemperedCfgs(\dCfgs)
  \right)
  = 1.
\end{equation}
We will call elements of $\wellTemperedCfgs(\dCfgs)$ \term{well-tempered configurations}.


\subsubsection{A weight for tempered configurations}
It will be convenient to compare two different configurations, especially in the dynamical setting we will see in the next section. To this end, we introduce (following Royer (\cite{Roy99})  the following weight:
\begin{equation}
  \label{eq=definitionAlpha}
  \alpha(i) = \sum \frac{p^{\convol n}}{\sigma'^n} (i),
\end{equation}
where $p^{\convol n}$ is the convolution product $p\convol p \cdots \convol p$,  $\sigma'$ satisfies $\sigma'<\inf V_1''$, and we recall that $p(i) = \abs{J(i)}$.
\begin{prpstn}
  The weight $\alpha$ decays exponentially: 
  \begin{equation}
    \label{eq=alphaDecay}
    \exists c_\alpha, d_\alpha, \quad
    \alpha(i) \leq c_\alpha \exp\left( - d_\alpha \abs{i} \right).
  \end{equation}

  Moreover, it satisfies the following~:
  \[
  \alpha \convol p \leq \sigma \alpha
  \]
\end{prpstn}
The proof is easy, and we refer to \cite{Zit06} for details.

The exponential decay of $\alpha$ shows that the tempered configurations $\spins$ have a finite $\lDeuxAlpha$-norm:
\[ \sum \alpha(i) \spin{i}^2  < \infty.\]

\subsubsection{The Glauber--Langevin dynamics}
It may be shown (\cf{} \cite{Roy99}), using standard tools, that in a finite volume $L$, the following SDE in $\xR^L$ has a strong solution~:
\begin{equation}
  \label{eq=SDE}
  dX_s = dB_s - \nabla U_{L,\spinsZ} (X_s) ds,
\end{equation}
where $B_s$ is a standard $\card(L)$-dimensional brownian motion. 

Using the $\lDeuxAlpha$ norm (where $\alpha$ is defined by \eqref{eq=definitionAlpha})
and  Gronwall-like arguments,  it is possible to compare the processes in different boxes, or starting from different points. 
We denote by $X^{L,\spinsZ,\spins}_t$ the process starting from $\spins$, in the box $L$ with boundary condition $\spinsZ$.
\begin{prpstn}
  For every set $L$, let $\beta_L(j) = \sum_i \ind{i\notin L} \alpha(i)\alpha(j-i)$.
Then,  for all $L\subset M$, and every tempered configuration $\spins$, 
  \begin{align}
    \label{eq=entreDeuxBoites}
    \alpha(0) \esp\left[ 
      \sum_i \alpha(i) \sup_{[0,t]} \left(X^{L,0,\spins} - X^{M,0,\spins}\right)^2 
    \right]
    & \leq e^{k't} \left( \norme{l^2(\beta_L)}{\spins}^2 
    + c\abs{\alpha}\sum_{j\in M\setminus L}\alpha(j)
      \right)
  \end{align}
  Moreover, for every  tempered $\spins$, $\norme{l^2(\beta_{L_n})}{\spins} \to 0.$
\end{prpstn}
The proof, inspired by \cite{Roy99}, may be found in \cite{Zit06} (lemmas 36 and 38).
These comparisons enable us to build an infinite volume dynamics. Moreover, we may let 
$M$ go to $\mathbb{Z}^d$ in \eqref{eq=entreDeuxBoites}
to get the following:
\begin{prpstn}
  Let $X^{0,\spins}$ be the infinite-volume process starting from $\spins$. Then:
  \begin{equation}
   \esp\left[ 
    \sum_i \alpha(i) \sup_{[0,t]} \left(X^{L,0,\spins} - X^{0,\spins}\right)^2 
    \right]
     \leq e^{k't} \left( \norme{l^2(\beta_L)}{\spins}^2 
    + c\abs{\alpha}\sum_{j\notin \setminus L}\alpha(j)
    \right)
  \end{equation}
\end{prpstn}
This gives an explicit estimate of the error made when we approximate the infinite volume dynamics by the finite volume one. This estimate can be made even more explicit if we use the decay properties of $\alpha$ (equation \eqref{eq=alphaDecay}, \cf{} lemma 38 and prop. 39 of \cite{Zit06}).
\begin{prpstn}\label{eq=FSP}
  There is a $d_\alpha$ such that:
  \begin{equation}
    \label{eq=entreFiniEtInfini}
    \forall \spins\in\temperedCfgs, \exists c_\spins, 
    \forall n, \quad
   \esp\left[ 
      \sum_i \alpha(i) \sup_{[0,t]} \left(X^{L,0,\spins} - X^{0,\spins}\right)^2 
    \right]
     \leq c_\spins \exp\left( k't - d_\alpha n \right)
  \end{equation}
\end{prpstn}

\subsection{Polynomial bounds on the entropy and related quantities}
We will need bounds on some entropy-related quantities in finite time. 
\begin{prpstn}
  \label{dfn=notationsEntropies}
  Let $\spins$ be a well-tempered configuration, and consider the processes $X^n\egalDef X^{L_n,0,\spins}_t$.
  We define the following notations:
  \begin{itemize}
    \item $h_t^n$ is the density of the law of $X^n_t$ with respect to the equilibrium measure 
      $\mu_{L_n, 0}$;
    \item $H_t^n$ is the entropy of this law ($H_t^n = \Ent_{\mu_{L_n}}(h_t^n)$);
    \item $H_{p,t}^n$ is given by
      \[ H_{p,t}^n = \int h_t^n \log_+^p( h_t^n) d\mu_n,\]
      where $\log_+$ is the positive part of the logarithm.
  \end{itemize}
  Then there exists a polynomial $Q$ (depending on $\spins$) such that, for all $p\geq 1$, and all $t\geq 1$, 
  \begin{equation}
    \label{eq=majorationPEntropie}
      H_{p,t}(n) \leq Q(n) p.
    \end{equation}
  In particular, since $H_t^n\leq H_{1,t}^n$, the entropy is polynomially bounded.

  Moreover, the degree of $Q$ does not depend on $\spins$ (as long as $\spins$ is well-tempered).
\end{prpstn}
This is a refinement of a result by Royer \cite{Roy99} (which deals only with the entropy, and does not 
precise the degree of $Q$, which will be needed later). The proof uses Girsanov's theorem to get an explicit expression of $h_t^n$, which is then estimated directly. The details may be found in \cite{Zit06}.

\section{From logarithmic Sobolev inequalities to uniqueness}
\label{sec=LS}
\subsection{Functional inequalities}
Let us start by recalling a few definitions.
\begin{dfntn}
  The measure $\mu$ on $\mathbb{R}^L$ satisfies a \term{logarithmic Sobolev inequality} with constant $C$ if 
  \begin{equation}
    \label{eq=dfnLS} \Ent_\mu(f^2) \leq C \int \sum_{i\in L} \abs{\nabla_i f}^2 d\mu.
  \end{equation}
  for every function $f$ such that both sides make sense. 
  
  It satisfies a \term{Poincaré inequality} with constant $C$ if 
  \begin{equation}
    \label{eq=dfnPoincare}
    \Var_\mu(f) \leq C  \int \sum_{i\in L} \abs{\nabla_i f}^2 d\mu.
  \end{equation}
  with the same restriction. 

  Finally, for $a\in(0,1)$, $\mu$ satisfies a \term{generalized Beckner inequality} \GBIa with constant $C$, if for any $f$, 
  \begin{equation}
    \label{eq=dfnBeckner}
    \sup_{p\in ]1,2[} \frac{ \int f^2 d\mu - (\int f^p d\mu) ^{2/p} } { (2-p)^a} \leq C_a\int \abs{\nabla f}^2 d\mu.
  \end{equation}
\end{dfntn}
The first two inequalities are well known, the third one was introduced in this form  by R.~Latala and K.~Oleszkiewicz in \cite{LO00}. It is known (\cf \cite{LO00,BCR05b}) that we recover Poincaré, resp. log-Sobolev,  by letting $a$ go to zero, resp. $1$, in the definition of the Beckner inequality. It is also known that Beckner inequalities may be compared~: \GBIa implies \GBIaPrime, whenever $a>a'$.

We will prove the uniqueness starting from hypotheses on the finite volume Gibbs measures, expressed in terms of functional inequalities. More precisely, we fix a boundary condition (for simplicity, we choose the $0$ boundary condition, however the same results should hold if we replace $0$ by a (fixed) tempered configuration $\spinsZ$), and make assumptions on the measures
$\mu_n = \mu_{L_n, 0}$.
\begin{assmptn}\label{assm=LS}
  $\mu_n$ satisfies a logarithmic Sobolev inequality, with constant $C_n$, where~:
  \[ C_n \leq C \frac{n}{\log(n)}, \]
  and $C$ is smaller than some explicit value (\cf{} \eqref{eq=Ccritique}).
\end{assmptn}
\begin{assmptn}\label{assm=Beckner}
  $\mu_n$ satisfies a Beckner($a$) inequality, with constant $C$, where $a$ and $C$ do not depend on $n$. Moreover, $a>a_{min}$, where $a_{min}$ only depends on the potential and the lattice dimension (\cf{} \eqref{eq=aEtD} for its explicit value).
\end{assmptn}
The main theorem  is the following.
\begin{thrm}
  \label{thm=uniciteParBeckner}
  If either Assumption \ref{assm=LS} or Assumption \ref{assm=Beckner} holds, there is only one tempered Gibbs measure in 
  infinite volume.
\end{thrm}

\subsection{Uniqueness from logarithmic Sobolev}
In this section, we prove theorem \ref{thm=uniciteParBeckner} under Assumption \ref{assm=LS}.

The main argument is the following. Let $P_t^L$ be the semi-group defined by the SDE \eqref{eq=SDE}
in the finite box $L$, with boundary condition $0$, and $P_t$ be the infinite-dimensional semi-group. For every $f$ (in a class to be precised later), we can decompose $P_t f$ in the following way:
\begin{equation}
  \label{eq=TheDecomposition}
  P_t f = \left(P_t f - P_t^L f\right)
 + \left( P_t^L f - \mu_L f\right)
 + \mu_L f.
 \end{equation}
 The first term may be controlled thanks to equation \eqref{eq=entreFiniEtInfini}.
 To get a good bound, we see that the diameter of $L$ should be of the order of $t$, to compensate the $\exp(kt)$. 

 More precisely, 
let us fix a $\rho$ (a ratio between $n$ and $t$) such that $\rho> k' >d'_\alpha$, and define $n(t) = \lfloor{\rho t}\rfloor +1$. By design, $n(t)$ satisifies:
\begin{equation}
  \label{eq=rapportEntreNetT}
  n(t) \in [\rho t, \rho t +1 ]
\end{equation}

This ensures
\[ k' t - d'_\alpha n(t) \leq (k' - \rho d'_\alpha) t,\]
where $k' - \rho d'_\alpha$ is a negative constant.
Hence when $t$ goes to $\infty$, $\esp\left[ \alphaNorm{\xBoite - \xInf}^2 \right]$ goes to zero. Plugging this back into \eqref{eq=entreFiniEtInfini} yields:
\[
\forall \spins\in\temperedCfgs,\quad
\abs{P_t f(\spins) - P_t^{L_{n(t)}}f(\spins)} \xrightarrow[t\to\infty]{} 0.
\]

\bigskip

The second term of \eqref{eq=TheDecomposition} depends on the convergence of a finite-dimensional diffusion to its equilibrium measure. This is where our functional inequalities come into play. Indeed, thanks to Pinsker's inequality and the exponential decrease of the entropy, 
\begin{align}
  \abs{ P_t^{n(t)} f - \mu_n f}^2 
  &\leq \osc^2(f) \norme{vt}{\loi{X_t^{L_{n(t)},0,\spins}} - \mu_{n(t)}}^2 \notag\\
  &\leq 2 \osc^2(f) 
    \relEnt{\loi{X_t^{L_{n(t)},0,\spins}}}{\mu_{n(t)}} 
         & &
  \notag \\
  &\leq 2 \osc^2(f) \exp\left( -2\frac{t-1}{\cNt}\right)
  \relEnt{\loi{X_1^{L_{n(t)}, 0,\spins}}}{\mu_{n(t)}}
     & & 
   \notag \\
\label{eq=convergenceBoiteFinie}
  &\leq 2 c_\spins \osc^2( f) \exp\left(-2\frac{t-1}{\cNt}\right) (1+n(t))^{d+d_\spins d_V}.
     & & \text{(prop. \ref{dfn=notationsEntropies})}.
\end{align}
\begin{rmrk}
  \label{rmq=casUniforme}
  Note that if we suppose (following Royer) a \emph{uniform} logarithmic Sobolev inequality, the proof is easily concluded: since $n$ is of the order of $t$, the power of $n$ is a power of $t$, and the exponential term ensures the convergence to zero. 
\end{rmrk}
Recall that tempered measures charge only well-tempered configurations. If we consider the left-hand side only for such configurations, we may replace $d_\spins$ by $\dCfgs$ on the right-hand side. 

Since by hypothesis, $C_{LS}(L_n) \leq C \frac{n}{\log(n)}$, and since $n(t)\leq \rho t + 1$,
\begin{align*}
 \exp\left( -2 \frac{t-1}{\cNt}\right) 
 &\leq \exp\left( -2 \frac{(t-1)\log(n(t))}{Cn(t)}\right)\\
 &\leq \exp\left( -2 \frac{(t-1)}{C (\rho t +1) } \log(n(t)) \right)
\end{align*}
For all $C'>C\rho$, and all $t$ large enough, $(t-1)/(C (\rho t + 1)) > 1/C'$, therefore
\begin{align*}
 \exp\left( -2 \frac{t-1}{\cNt}\right) 
 &\leq \exp\left( - \frac{2}{C'} \log(n(t)) \right)\\
 &\leq n(t)^{-2/C'}
\end{align*}
Coming back to  \eqref{eq=convergenceBoiteFinie}, we obtain 
\begin{align*}
  \abs{ P_t^{n(t)} f - \mu_n f}^2 
  &\leq c'_\spins \osc^2(f) n(t)^{-2/C'} (1+n(t))^{d + \dCfgs d_V}
\end{align*}
The r.h.s. converges as soon as         $d + \dCfgs d_V<2/C'$, \ie{}:
\[ C' < \frac{2}{d+\dCfgs d_V}.\]
This is possible if 
\begin{equation}
  \label{eq=Ccritique}
  C<\frac{2}{\rho(d+\dCfgs d_V)}.
\end{equation}

Under this condition, we have shown:
\[
\forall \spins \in \temperedCfgs,\quad
 \abs{P_t^{n(t)} f(\spins) - \mu_n f} \xrightarrow[t\to\infty]{} 0.
\]

Once we have chosen a scale $n(t)$ that guarantees convergence for the first two terms of \eqref{eq=TheDecomposition}, the last one may be dealt with thanks to a compacity argument (\cite{Roy99} p.~72)
 \[ \exists (t_k), t_k \to\infty, \exists \mu\in\temperedDLRStates, 
 \mu_{L_{n(t_k),0}} \xrightarrow{k\to\infty} \mu.
   \]

  Along this  particular sequence $t_k$ of times, 
 \[
   \forall \spins, \spins \text{ tempered}, \forall f, 
   \quad P_{t_k}f(\spins) \xrightarrow{k\to\infty} \mu f.
 \]
 Let then $\nu$ be another tempered Gibbs measure. It is known (\cf{} \cite{Roy99}, Theorem 4.2.13)
  that $\nu$ is necessarily invariant w.r.t the semi-group $P_t$. 
 Then 
 \[
   \nu(f) = \nu (P_{t_k} f).
 \] 
 $P_{t_k}f$ converges pointwise to $\mu f$. Letting $k$ go to infinity, since $f$ is bounded, we have by dominated convergence:
 \[
   \nu(f) = \nu( \mu f) = \mu f.
 \]
 Since this is true for $f$ in a sufficiently large class of functions, $\nu = \mu$ and the tempered Gibbs measure is unique: theorem \ref{thm=uniciteParBeckner} follows from Assumption \ref{assm=LS}. 


\section{Beyond logarithmic Sobolev inequalities}

\subsection{Strong enough Beckner inequalities imply uniqueness}
We now prove uniqueness under assumption \ref{assm=Beckner}.
The compacity argument and the comparison between finite and infinite volume still hold~; 
the only thing to check is that the assumption is strong enough to guarantee~:
\[
  P_t^L f - \mu_L f \xrightarrow{t\to\infty} 0.
\]
Once more, we use Pinsker's inequality to bound this difference by an entropy. 
This entropy does not decay exponentially fast (since we do not suppose log-Sobolev inequalities anymore), but we are able to show that 
it converges nonetheless. 

The argument is adapted from \cite{CGG05}, where the following result is proved.
\begin{thrm}
  [\cite{CGG05}, Th. 5.5]
  \label{thm=decroissanceCGG}
  Let $\mu$ be a probability measure on  $\xR^n$, absolutely continuous w.r.t. Lebesgue measure, 
  and satisfying a \GBIa{} inequality.

  Let $\nu$ (an initial law) be such that:
  \begin{equation}
    \label{eq=hypotheseCGG}
    \exists t_0,C , \forall t\geq t_0, \forall p\geq 1, \quad
  \left( \int P_t\nu \log_+^p (P_t \nu) \right)^{1/p} \leq Cp.
\end{equation}
Then the entropy starting from $\nu$ decays sub exponentially along $P_t$:
  \begin{equation}
    \label{eq=decroissanceCGG}
  \forall a'<a, \exists s,t_0,  \forall t\geq t_0, \quad
  \Ent_\mu(P_{s+t}\nu) \leq \exp\left( 1 - t^{1/(2-a')}\right).
  \end{equation}
\end{thrm}
Note that our parameter $a$ is linked to the $\alpha$ appearing in \cite{CGG05} by $a = (2\alpha - 2) /\alpha$ (\cf{} example 4.3 of \cite{CGG05}).

Since this theorem entails a fast convergence (faster than polynomial), it is natural to expect that it should be enough for our purposes. Unfortunately, this results holds for large (and unspecified) $t$, and we need to use it for a relatively small $t$ (of the order of the diameter of the box $L$).

Therefore, we will use ideas of \cite{CGG05} to prove a similar result with explicit constants. The 
preliminary estimates we need were already cited in the previous section (\cf{} \eqref{eq=majorationPEntropie})
. We prove the result in two steps: first we bound the entropy for small times, then we iterate the estimate.

\subsubsection{First entropy estimate}
We follow the idea of the fifth section of \cite{CGG05}, (\emph{Convergence to equilibrium for diffusion processes}).

It is well-known that logarithmic Sobolev inequalities imply exponential convergence of the entropy.
When we only have a Beckner inequality, we may still prove exponential convergence, but only for bounded functions. 
\begin{lmm}[\cite{CGG05}, example 4.3]
  If $\mu$ satisfies \GBIa, there exists $C'_a$ such that, for any \emph{bounded} probability density $h$:
  \begin{equation}
    \label{eq=decroissanceEntropieBornee}
    \Ent_\mu(P_t h) \leq \Ent_\mu(h) \times \exp\left(
      - \frac{t}{C'_a \left( 1 + \log^{1-a}\left( \ninf{h}\right)\right)}
      \right).
    \end{equation}
\end{lmm}
Cattiaux, Gentil and Guillin have shown that this implies decay estimates for all functions, but the decay is not exponential. 

The idea of their proof is to decompose a function $h$ in a bounded part $h\ind{h\leq K}$ and a remainder $h\ind{h>k}$, and then choose an appropriate $K$.

Using this method, we prove the following:
\begin{prpstn}
  \label{lmm=controleHst}
  If $\mu_n$ satisfies a \GBIa{} inequality,  uniformly in  $n$, 
  then for all $a<a_0$, 
  there exists a polynomial $Q = Q_{a,\spins}$, whose degree depends only on $V$ and the dimension $d$, and a number $t_0(a)$, such that
  \begin{equation}
    \label{eq=controleHst}
    \forall s\geq 1, \forall t\geq t_0(a), \quad
  H_{s+t}^n \leq \frac{1}{c_{t,n}} \phi(H_s),
  \end{equation}
  where $\phi(x) = x \left( 1 + \log_+(1/x) \right)$,
  and $c_{t,n} = t^{1/(1-a)} / Q(n)$.
\end{prpstn}
We will need the following lemma, which we quote without proof.
\begin{lmm}
  [\cite{CGG05}, lemma 5.3]
  \label{lmm=controleEntropieTronquee}
  Let $h$ be a probability density w.r.t. $\mu$. If there exists  $c>0$ such that the $p$-entropy is bounded:
  \[ \forall p>1, \quad H_{p,t} \leq cp,\]
  and if $K$ satisfies
  \[ K\geq e^2, \qquad \log(K) \geq 2e\times\Ent_\mu( h ),\]
  then
  \[
  \Ent_\mu(h\ind{h>K}) \leq (ec + 2) \frac{\Ent_\mu(h)}{\log(K)} \log\left(
    \frac{\log(K)}{\Ent_\mu(h)} \right).
  \]
\end{lmm}

We will also need bounds on the entropy of bounded functions.
\begin{lmm}
  \label{lmm=entropieFonctionsBornees}
  Let $\mu$ be a measure that satisfies \GBIa. There exists $C'_a$ such that, if $h$ is a bounded 
  probability density, $H$ is the entropy of $h$, and $K$ satisfies
  \[ K\geq e^2, \qquad \log(K) \geq 4 H,\]
  then
  \[
  \Ent_\mu\left( P_t (h \ind{h\leq K} ) \right)
  \leq H\times  \exp\left( \frac{t}{ C'_a \log^{1-a}\left( K \right) } \right).
  \]
\end{lmm}
\begin{proof}
  This lemma follows from equation
  \eqref{eq=decroissanceEntropieBornee}. 
  In order to see this, we would like to normalize
  $h\ind{h\leq K}$ so that it becomes a probability density, 
  and apply the previous lemma. 
  This can be done if $\int h \ind{h\leq K} \neq 0$.
  Lemma 3.4 from \cite{CG06}, shows that, for $K\geq e^2$~:
  \[
  \int h \ind{h>K} \leq \frac{2H}{\log K}.
  \]
 Since we assume  $\log(K)\geq 4H$,
  \[
  \int h\ind{h>K} \leq \frac{1}{2},
  \]
  and since $\int h = 1$, this entails:
  \begin{equation}
    \label{eq=minorationIntegraleH}
    \int h \ind{h\leq K}  = 1 - \int h\ind{h>K} \geq 1/2.
  \end{equation}
  Let us denote by $\tilde{h}$ the renormalized version of $h\ind{h\leq K}$. It is a bounded probability density, and we may apply the bound \eqref{eq=decroissanceEntropieBornee}~:
  \[
  \Ent_\mu ( P_t \tilde{h} ) 
  \leq \exp \left( - \frac{t}{C'_a \left(
      1 + \log^{1-a}\left( \ninf{\tilde{h}} \right) \right)} 
    \right) \Ent_\mu(\tilde{h}).
  \]
  We multiply both sides by $\int h\ind{h\leq K} d\mu$, and put these factors in the entropies (by homogeneity).
  \begin{equation}
    \label{eq=decroissancePtHK}
  \Ent_\mu (P_t (h\ind{h\leq K}))
  \leq \exp \left( - \frac{t}{C'_a \left(
      1 + \log^{1-a}\left( \ninf{\tilde{h}} \right) \right)} 
    \right) \Ent_\mu(h\ind{h\leq K}).
  \end{equation}
  Finally, we control the sup norm of $\tilde{h}$:
  \begin{align*}
    \tilde{h}  &= \frac{ h \ind{ h\leq K} } { \int h \ind{h\leq K} d\mu } \\
  & \leq \frac {K}{ 1/2},
  \end{align*}
  where we reused the bound  \eqref{eq=minorationIntegraleH} on the integral.
  Since $K\geq e$, the denominator of \eqref{eq=decroissancePtHK} is bounded above:
  \begin{align*}
  C'_a
    &\left( 1+ \log^{1-a}(\ninf{\tilde{h}})\right)  \\
    &\leq C'_a \left( \log^{1-a}(K) + \log^{1-a}(2K) \right)  \\
    &\leq C''_a \log^{1-a}(K).
  \end{align*}
  This proves the lemma.
\end{proof}

We now proceed to the proof of the proposition  \ref{lmm=controleHst}.
\begin{proof}
  By definition, $H_t = \Ent_{\mu_{L_n}}(h_t)$. We consider the time $s+t$, and truncate $h_s$: for all $K$,
  \[ h_s = h_s \ind{h_s \leq K} + h_s \ind{h_s >K}.\]
  For any positive functions $(f,g)$, $\Ent(f+g) \leq \Ent(f) + \Ent(g)$ --- this follows easily  from the variational formula for the entropy: $\Ent_\mu(f) = \sup\left\{ \int fh, \int e^h d\mu = 1 \right\}$. Therefore,
  \begin{align*}
    \forall s,t, \forall K,\quad
    H_{t+s} =  \Ent (P_t h_s) 
    &\leq \Ent( P_t (h_s \ind{h_s \leq K})) + \Ent(P_t (h_s \ind{h_s > K})) \\
    &\leq \Ent( P_t (h_s \ind{h_s \leq K})) + \Ent(    h_s \ind{h_s > K}),
  \end{align*}
  since the entropy decreases along $P_t$.
  Suppose that $K$ satisfies:
  \begin{equation}
    \label{eq=conditionsK}
  \begin{cases} K\geq e^e, \\
     \log(K)\geq 2e  H_s.
   \end{cases}
 \end{equation}
   We now apply lemma \ref{lmm=entropieFonctionsBornees} to the first term. For the second term, \eqref{eq=majorationPEntropie} shows that the hypotheses of lemma \ref{lmm=controleEntropieTronquee} are fulfilled. 
   If $K$ satisfies both hypotheses, we get:
  \begin{equation}
    \label{eq=majorationHts1}
  H_{t+s} 
  \leq \exp\left(
    - \frac{t}{C_{a_0}\log(K)^{1-a_0}} \right)H_s + Q(n) \frac{H_s}{\log(K)} \log\left( \frac{\log(K)}{H_s}
    \right).
\end{equation}
  We now define $K$ to be the unique solution on $(e^e,\infty)$ of the following equation: 
  \begin{equation}
    \label{eq=choixDeK}
    \log(K)  = \left( \frac{ t }{C_{a_0} \log\log K} \right)^{1/(1 -a_0)},
  \end{equation}
  This $K$ (which depends on $t$) is well defined, because
   $K\mapsto \log(K) \log\log(K)^{1/(1-a_0)}$ is bijective from $]e^e, \infty[$ onto $]0, \infty[$.
  Assume for the time being that $K$ satifies the second condition of \eqref{eq=conditionsK}.
  The inequality \eqref{eq=majorationHts1} becomes
  \[
  H_{t+s} \leq 
    \frac{1}{\log(K)} H_s
    + Q(n) \frac{H_s}{\log(K)} \log\left( \frac{\log(K)}{H_s}\right). 
  \]
  Let us work a little bit to get a simpler bound. Since $K\geq e^e$, $\log \log K\geq 1$, and $Q(n)$ may always be taken larger than $1$. This yields:
  \begin{align}
    \notag
    H_{t+s} &\leq \frac{\log\log(K)}{\log(K)} H_s + \frac{Q(n)H_s}{\log(K)} \left( \log\log(K) + \log_+(1/H_s) \right) \\
    \notag
    &\leq         \frac{\log\log(K)}{\log(K)} H_s + \frac{Q(n)H_s}{\log(K)} \log\log(K) \left( 1 + \log_+(1/H_s)\right) \\
    \label{eq=majorationHts2}
    &\leq \left( \frac{\log\log(K)}{\log(K)}\right) Q(n) H_s \left( 2 + \log_+\left( 1/H_s\right) \right).
  \end{align}
  Our choice of $K$ ensures that there exists a constant $c_a$ such that:
  \[ \log\log(K)/\log(K) \leq \frac{c_a}{ t^{1/(1-a)}}, \]
  for any $t$ larger than a $t_0(a)$.

  Insert this into equation \eqref{eq=majorationHts2}, and define $Q_a(n) = 2c_a Q(n)$. The bound becomes
  \[
   H_{t+s} \leq \frac{ Q_a(n) } {t^{1/(1-a)}} H_s \left( 1 + \log_+(1/H_s) \right),
   \]
   which is exactly \eqref{eq=controleHst}.

   Let us go back to the case where $K$ (defined as the solution of \eqref{eq=choixDeK}) does not satisfy \eqref{eq=conditionsK}. Since by definition $K\geq e^e$, we need only consider the case where $\log(K) \leq 2eH_s$. 
   We know that there exists a polynomial $Q'_{a,\spins}$ such that $H_s^n \leq Q'_{a,\spins}(n)$, for all $s \geq 1$. In this case, 
   \[ \log(K) \leq Q''_{a,\spins}(n).\]
   In other words,
   \[ 1 \leq \frac{Q''_{a,\spins}(n)}{\log(K)}.\]
   Since $\log\log(K) \geq 1$, it follows that
   \[ 1 \leq \frac{\log\log K}{\log K} Q''_{a,\spins}(n).\]
   Finally, the entropy $H$ decreases along the semi-group. For every $s\geq 1$, and $t\geq t_0(a)$, we have:
   \[
   H_{t+s} \leq H_s 
   \leq \left(
       \frac{ \log \log K}{\log K} 
     \right) Q''_{a,\spins}(n) H_s (2 + \log_+(1/H_s)).
   \]
   This shows that \eqref{eq=majorationHts2} still holds, and the end of the proof is the same.
\end{proof}

\subsection{Iteration of the estimate} 
We are now in a position to prove theorem \ref{thm=uniciteParBeckner}, under Assumption \ref{assm=Beckner}.

The previous estimate \eqref{eq=controleHst} is useful if $c_{t,n}$ is greater than $1$. Let $D$ be the degree of $Q$ (it does not depend on $\spins$ nor on $a$). We will assume:
\begin{equation}
  \label{eq=aEtD}
  a_0 > a_{min} = \frac{D-1}{D}.
\end{equation}
Note that we may choose $a>a_{min}$ in lemma \ref{lmm=controleHst}.
\begin{lmm}
  \label{lmm=presqueFini}
  The following properties hold.
  \begin{itemize}
    \item There exists $t_0(n)$ such that, for all $t> t_0(n)$, $c_{t,n} > 1$~;
    \item There exists $u_0(n)$ such that, for all $u>u_0$, $u$ may be written as $t(c_{t,n})^2$, with  $t>t_0(n)$~;
    \item The quantity $u_0(n)$ is relatively small:
      \begin{equation}
	\label{eq=u0estPetit}
	u_0(n) = o(n).
      \end{equation}
 \end{itemize}
 
 Moreover, for all $u>u_0(n)$,for all $s\geq 1$, 
 \begin{equation}
   \label{eq=decroissanceHu}
 H^n_{s+u} \leq (e+H_s) \exp\left( -\frac{ u^{1/(3-a)}} {Q(n)^{( 1 - a)/(3-a)}} \right).
 \end{equation}
\end{lmm}
This lemma implies theorem \ref{thm=uniciteParBeckner}.

Indeed, we only need to show that the entropy at time $t$ in the box $L_{n(t)}$ converges to zero. We apply the lemma with $s=1, u=t, n = n(t)$ (this is possible thanks to \eqref{eq=u0estPetit}).
Since $Q(n)$ is (by definition) of degree $D$, it is bounded above by $n^D$ (up to a constant), and there is a $c$ such that:
\[
H^{n(t)}_{1+t} 
\leq \left(e+ H^{n(t)}_1\right) \exp \left( -c \left(\frac{n(t)}{n(t)^{D(1-a)}}\right)^{1/(3-a)} \right).
\]
Since $H^{n(t)}_1$ grows polynomially in $n$ (this is the result of theorem \ref{dfn=notationsEntropies}) and $t$ is of the order of $n$, it suffices to show that the power of $n$ in the exponential is positive, and the whole quantity will go to zero. This power is:
\[ \frac{1}{3-a}\left( 1-D+aD \right),\]
which is indeed positive, because $a>a_{min}$ (defined by \eqref{eq=aEtD}).

This shows that the entropy at time $t$ in the box $L_{n(t)}$ goes to zero. As was already said before, the other parts of the proof require no change, therefore theorem \ref{thm=uniciteParBeckner} will be proved as soon as we show lemma \ref{lmm=presqueFini}.

\begin{proof}[\proofname{} of lemma \ref{lmm=presqueFini}]
  Let us begin by showing the existence of $t_0$, $u_0$.

  Recall that $c_{t,n} = t^{1/(1-a)}/Q(n)$, and that the degree of $Q$ is $D$. Let us choose an $a'$ such that
  \[ a_{min} < a' < a < a_0, \]
  If we define $t_0 = cn^{D(1-a')}$ for some constant $c$, 
  \[ c_{t,n} \geq \frac{t^{1/(1-a)}} {n^D} >1, \]
  for $t\geq t_0$.

  One then defines $u=u(t,n) = t(c_{t,n})^2$. This increases with $t$, and one may choose
  \[
  u_0(n) = u( t_0(n),n)  = c n^{D(1-a')(3-a)/(1-a)} / Q(n)^2. 
  \]

  We would like $u_0$ to be small w.r.t. $n$ (we do not want to wait for a period longer than the diameter of the box). The previous choice ensures:
  \[ 
  u_0(n) \sim c n^{D(1-a')(3-a)/(1-a) - 2D}.
  \]
  This $u_0$ is negligible compared with $n$
  when
  \[ D\frac{(1-a')(3-a)}{1-a} - 2D <1.\]
  This is satisfied for $a'=a$ (because $D(3-a) - 2D = D(1-a) <1$, since $a>a_{min}$). By continuity, this still holds for some $a'<a$.

  \bigskip

  Let us now prove \eqref{eq=decroissanceHu}. The idea is to iterate the estimate given by lemma \ref{lmm=controleHst}. To do this, fix $t>t_0(n)$, and define the sequence $(u_k)$ by $u_k = H_{s+kt}$.
  To control $u_k$, we compare it to $v_k$ defined recursively by:
  \[
  \begin{cases}
    v_0 = u_0,\\
    v_{k+1} = f(v_k),
  \end{cases}
  \]
  where $f(x)=\frac{1}{c_{t,n}} \phi(x)$ (\cf{} lemma    \ref{lmm=controleHst}).
  Since $f$ is increasing, and
  \[ u_{k+1} \leq f(u_k) \]
  (this follows from equation \eqref{eq=controleHst}, applied with $s=s+tk$, and $t=t$), 
  it is easily seen by induction that $u_k \leq v_k$.
  
  Now $v_k$ is easily studied by standard methods: the condition $c_{t,n}>1$ ensures that $f$ has only one stable stationary point, $x_e = \exp( 1 - c_{t,n})$, and that $v_k$ converges to this point.
  If we start from a point to the left of $x_e$, $v_k$ is always bounded by $x_e$.
  On the right of $x_e$, $f$ is a  $\left( 1 - \frac{1}{c_{t,n}}\right)$-contraction. Therefore, for all $k$,
  \begin{equation}
    \label{eq=iterationPourV}
    v_k \leq x_e + \left(1 - \frac{1}{c_{t,n}}\right)^k \left(v_0 - x_e\right)_+.
  \end{equation}
  The explicit value of $x_e$, and the bound      $(1 - 1/c)^k \leq \exp ( - k/c)$ show that:
  \[ 
  \forall k, \quad v_k \leq \exp(1 - c_{t,n}) + v_0 \exp ( -k/c_{t,n}).
  \]
  Let us now look at the entropy $H$ at time   $s+u$.  If $u$ can be written as $u=kt$ with a  $t>t_0(n)$, the previous iterated bound reads:
  \[ 
  H_{s+u} \leq \exp( 1 - c_{t,n}) + v_0 \exp ( -k/c_{t,n}).
  \] 
  For any $u>u_0(n)$, let us choose $t$ such that $t(c_{t,n})^2 = u$, and $k = c_{t,n}^2$ (more precisely, $k$ is the nearest integer). By definition of $u_0$, $t$ is larger than $t_0$. Now $t$ and $c_{t,n}$ may be rewritten as functions of $u$:
  \begin{align*}
    t c_{t,n}^2 = u  \text{, therefore } t &= \left(uQ(n)^2\right)^{(1-a)/(3-a)}, \\
    c_{t,n} &= \frac{1}{Q(n)} (uQ(n)^2)^{1/(3-a)}. \\
  \end{align*}
  Since $u_0 = H_s$, we obtain:
\begin{align*}
  H_{s+u} &\leq (e + H_s) \exp\left( - c_{t,n}\right) \\
  & \leq (1+H_s)\exp\left( -\frac{ u^{1/(3-a)}} {Q(n)^{( 1-a)/(3-a)}} \right).
 \end{align*}
 This concludes the proof.
\end{proof}



\subsection{Are Beckner inequalities the right scale~?}
We show here that the scale of Beckner inequalities is arguably the ``right'' one for proving uniqueness. 
We only give a heuristic argument, using a toy model introduced by T.~Bodineau and F.~Martinelli in \cite{BM02}.

This paper studies the phase transition regime, and tries to find lower bounds on the growth of the constants, as the size increases.
The type of result they get is:
\begin{prpstn}
  In the phase transition regime, for the $+$ boundary condition, the LS constants (in  $[-n,n]^d$)  grow at least like $n^2$.
\end{prpstn}
This result is similar to our theorem:
If the proposition holds in our setting, and in the whole phase transition regime, then a sublinear growth of the LS constants must imply that no phase transition occurs.

 Their approach is however very different: they find a ``good'' test function for which the entropy is large whereas the energy stays small.  

In another section, the authors introduce a toy model, which is supposed to reproduce the main aspects of the dynamics for the (classical) Ising model, in the phase transition regime: namely, the dynamics of the disappearance of a big droplet of $-$ spins when the boundary condition is $+$. 

The model is a birth and death process on $\{0, n^d\}$.  with rates $b$ and $d$:
\begin{equation}
  \label{eq=definitionModele}
\begin{aligned}
  b(x)&= x^\alpha & &\text{if $x\geq 1$}, \\
  b(0)&= 1, &&\\
  d(x+1)&= x^\alpha \exp\left( (x+1)^\alpha - x^\alpha \right) & &\text{if $x\geq 2$}\\
  d(1)  &= e. &&
\end{aligned}
\end{equation}
We choose $\alpha=\frac{d-1}{d}$ and note that the process is reversible w.r.t. $\mu$ defined by $\mu(x) = \frac{1}{Z} \exp( -x^\alpha )$.

The authors of \cite{BM02} then proceed to study the Poincaré and log-Sobolev constants of this one-dimensional  by means of Muckenhoupt-like criteria, established in the discrete case by Miclo (\cite{Mic99}). 
In fact, similar results exist for \emph{any} Beckner inequality. We rephrase here a result from  \cite{BR03} (the discrete version of Theorem 13, justified in the remarks at the end of section 4 --- note that our $a$ is related to their $r$ by $a = 2(1 - 1/r)$).
\begin{prpstn}
  For any $i \in \xZ$, define the following quantities:
  \begin{align*}
    R_+(x) &= \sum_{y\geq x} \mu(y),  &  R_-(x) &= \sum_{y \leq x} \mu(y), \\
    S_+(i,x)&= \sum_{y=i+1}^x \frac{1}{\mu(y) b(y)}, &  S_-(i,x) &= \sum_{x}^{i-1} \frac{1}{\mu(y) b(y) } \\
    B_+(i) &= \sup_{x>i} S_+(i,x) R_+(i,x) \log^a \left(1 + \frac{1}{2R_+(i,x)} \right) &&\\
    B_-(i) &= \sup_{x<i} S_-(i,x) R_-(i,x) \log^a \left(1 + \frac{1}{2R_-(i,x)}\right).&&
  \end{align*}
  Finally, let $B = \inf_{i} (B_+(i) \wedge B_- (i))$.
  Then $\mu$ satisfies a \GBIa inequality if and only if $B$ is finite, and there exists a universal constant $k$ such that $\frac{1}{k} B \leq C_a \leq kB$.
\end{prpstn}

This can be used to find explicit bounds on the Beckner constants, thanks to the estimates (\cite{BM02}):
\begin{align*}
  \sum_{y\geq x} \mu(y) &\approx x^{1-\alpha} \exp( - x^\alpha), \\
 \sum_{y=i+1}^x \frac{1}{\mu(y) b(y) } &\approx x^{1-2\alpha} \exp (x^\alpha),
 \end{align*}
 where $X\approx Y$ means that there exists a $k$ (independant of $x,i,\alpha$) such that $\frac{1}k X \leq Y \leq kX$. This implies estimates on $B_+,B_-$ and $B$, \eg{}: 
 \begin{align*}
   B_+(i,x) &\approx x^{1-2\alpha} \exp(x^\alpha) x^{1-\alpha} \left(x^\alpha\right)^a \\
   &\approx x^{1-3\alpha + \alpha a}.
 \end{align*}
 Define $a_d$ to be the solution of $1 - 3\alpha + \alpha a = 0$. If $a>a_d$, $B \approx B_+(i,N^d) \approx N^{d(1 - 3\alpha + \alpha a)}$ and the Beckner constant blows up with $N$. If $a<a_d$, $B$, and therefore the Beckner constant, stays bounded with $N$. 

 Since $a_d = (3\alpha -1)/\alpha$ and $\alpha$ is defined in terms of a ``dimension'' $d$, we have shown the following
\begin{thrm}
  Consider the toy model defined by \eqref{eq=definitionModele}
  For each value of the ``dimension'' $d$, there exists an $a_d$ such that~:
  \begin{itemize}
    \item If $a>a_d$, the Beckner constant $C(a,N)$ grows like $N^{}$;
    \item If $a<a_d$, the Beckner constant $C(a,n)$ stays bounded in $N$.
  \end{itemize}
  Moreover, $a_d$ satisfies:
  \begin{itemize}
    \item If $d=1$ or $2$, $a_d<0$ so that all constants blow up in $N$;
    \item If $d=3$, $a_d=0$, the Poincaré constant stays bounded whereas all other constants blow up;
    \item If $d>3$, $a_d\in(0,1)$.
  \end{itemize}
\end{thrm}
In particular, this tells us that (if the toy model is an appropriate approximation of the true model), there may be parameters for which the phase transition occurs, but the Poincaré constant stays bounded. 
This leads us to believe that theorem \ref{thm=uniciteParBeckner}
should not be too far from optimality, and that it should not be possible to prove uniqueness if we only suppose a uniform Poincaré inequality.


\end{document}